\documentclass[12pt,reqno]{amsart}
\usepackage{amscd,amsmath,amsthm,amssymb}
\usepackage{color}
\usepackage{pstricks}
\usepackage{stmaryrd}
\usepackage{tikz}
\usepackage{url}

\usepackage{latexsym}
\usepackage{amsfonts,amsmath,mathtools}
\usepackage{graphics}
\usepackage{float}
\usepackage{enumitem}

\usepackage{booktabs} 
\usepackage{colortbl}
\usepackage{lipsum}

\newpsstyle{fatline}{linewidth=1.5pt}
\newpsstyle{fyp}{fillstyle=solid,fillcolor=verylight}
\definecolor{verylight}{gray}{0.97}
\definecolor{light}{gray}{0.9}
\definecolor{medium}{gray}{0.85}
\definecolor{dark}{gray}{0.6}

\def\NZQ{\mathbb}               

\def\QQ{{\NZQ Q}}
\def\ZZ{{\NZQ Z}}

\def\FF{{\NZQ F}}

\def\frk{\mathfrak}               

\def\mm{{\frk m}}

\def\G{{\mathcal G}}

\def\HS{\textup{HS}}
\def\pd{\textup{proj}\phantom{.}\!\textup{dim}}

\def\opn#1#2{\def#1{\operatorname{#2}}} 
\opn\chara{char} \opn\length{\ell} \opn\pd{pd} \opn\rk{rk}
\opn\projdim{proj\,dim} \opn\injdim{inj\,dim} \opn\rank{rank}
\opn\depth{depth} \opn\grade{grade} \opn\height{height}
\opn\embdim{emb\,dim} \opn\codim{codim}

\opn\Tr{Tr} \opn\bigrank{big\,rank}
\opn\superheight{superheight}\opn\lcm{lcm}
\opn\trdeg{tr\,deg}
	\opn\reg{reg} \opn\lreg{lreg} \opn\ini{in} \opn\lpd{lpd}
	\opn\size{size} \opn\sdepth{sdepth}
	\opn\link{link}\opn\fdepth{fdepth}\opn\lex{lex}
	\opn\tr{tr}
	\opn\type{type}
	\opn\gap{gap}
	\opn\diam{diam}
	\opn\Mod{Mod}
	\opn\div{div} \opn\Div{Div} \opn\cl{cl} \opn\Cl{Cl}
	\opn\Spec{Spec} \opn\Supp{Supp} \opn\supp{supp} \opn\Sing{Sing}
	\opn\Ass{Ass} \opn\Min{Min}\opn\Mon{Mon}
	\opn\Ann{Ann} \opn\Rad{Rad} \opn\Soc{Soc}
	\opn\Im{Im} \opn\Ker{Ker} \opn\Coker{Coker} \opn\Am{Am}
	\opn\Hom{Hom} \opn\Tor{Tor} \opn\Ext{Ext} \opn\End{End}
	\opn\Aut{Aut} \opn\id{id}
	
	\opn\nat{nat}
	\opn\pff{pf}
	\opn\Pf{Pf} \opn\GL{GL} \opn\SL{SL} \opn\mod{mod} \opn\ord{ord}
	\opn\Gin{Gin} \opn\Hilb{Hilb}\opn\sort{sort}
	\opn\PF{PF}\opn\Ap{Ap}
	\opn\dist{dist}
	\opn\aff{aff}
	\opn\relint{relint} \opn\st{st}
	\opn\lk{lk} \opn\cn{cn} \opn\core{core} \opn\vol{vol}  \opn\inp{inp} \opn\nilpot{nilpot}
	\opn\link{link} \opn\star{star}\opn\lex{lex}\opn\set{set}
	\opn\width{wd}
	\opn\Fr{F}
	\opn\QF{QF}
	\opn\G{G}
	\opn\type{type}\opn\res{res}
	\opn\conv{conv}
	\opn\sr{sr}
	\opn\gr{gr}
	
	\def\pot#1#2{#1[\kern-0.28ex[#2]\kern-0.28ex]}

	\opn\dirlim{\underrightarrow{\lim}}
	\opn\inivlim{\underleftarrow{\lim}}
	\def\Implies{\ifmmode\Longrightarrow \else
		\unskip${}\Longrightarrow{}$\ignorespaces\fi}
	\def\implies{\ifmmode\Rightarrow \else
		\unskip${}\Rightarrow{}$\ignorespaces\fi}
	\def\iff{\ifmmode\Longleftrightarrow \else
		\unskip${}\Longleftrightarrow{}$\ignorespaces\fi}

	\let\:=\colon
	\newtheorem{Theorem}{Theorem}[section]
	\newtheorem{Lemma}[Theorem]{Lemma}
	
	\newtheorem{Proposition}[Theorem]{Proposition}

	\newtheorem{Definition}[Theorem]{Definition}

	\let\epsilon\varepsilon
	\let\kappa=\varkappa
	\def\qed{\ifhmode\textqed\fi
		\ifmmode\ifinner\hfill\quad\qedsymbol\else\dispqed\fi\fi}
	\def\textqed{\unskip\nobreak\penalty50
		\hskip2em\hbox{}\nobreak\hfill\qedsymbol
		\parfillskip=0pt \finalhyphendemerits=0}
	\def\dispqed{\rlap{\qquad\qedsymbol}}
	
	\opn\dis{dis}
	\def\pnt{{\raise0.5mm\hbox{\large\bf.}}}
	
	\opn\Lex{Lex}
	\opn\Max{Max}
	\opn\Shad{Shad}
	\opn\astab{astab}

	\opn\v{v}
	\def\soc{\textup{soc}}
	\def\indeg{\textup{indeg}}
	
\begin{document}
	
\title{Homological Shift Ideals: Macaulay2 Package}	
\author{Antonino Ficarra}

\address{Antonino Ficarra, Department of mathematics and computer sciences, physics and earth sciences, University of Messina, Viale Ferdinando Stagno d'Alcontres 31, 98166 Messina, Italy}
\email{antficarra@unime.it}

\thanks{
}



\maketitle

\begin{abstract}
	We introduce the \textit{Macaulay2} package \texttt{HomologicalShiftIdeals}. It allows to compute the homological shift ideals of a monomial ideal, and to check the homological shift properties, including having linear resolution, having linear quotients, or being polymatroidal. The theory behind these concepts is explained and the main features of the package are presented.
\end{abstract}

\section{Introduction}

Let $S=K[x_1,\dots,x_n]$ be the standard graded polynomial ring with coefficients in a field $K$. Let $I$ be a monomial ideal of $S$. For a vector ${\bf a}=(a_1,\dots,a_n)\in\ZZ_{\ge0}^n$, we set ${\bf x^a}=x_1^{a_1}\cdots x_n^{a_n}$. Note that, as a $S$-module, $I$ is multigraded. Hence, the minimal free resolution of $I$, $\FF:\cdots\rightarrow F_i\rightarrow\cdots\rightarrow F_0\rightarrow I\rightarrow0$, is naturally multigraded. Thus, $F_i=\bigoplus_{\bf a}S(-{\bf a})^{\beta_{i,{\bf a}}(I)}$ for all $i$, where $\beta_{i,{\bf a}}(I)$ is a multigraded Betti number. The $i$th \textit{homological shift ideal} of $I$ is the monomial ideal defined as $$\HS_i(I)\ =\ ({\bf x^a}\ :\ \beta_{i,{\bf a}}(I)\ne0).$$

Homological shift ideals have been introduced in \cite{HMRZ021a}, and attracted the interest of many researchers \cite{Bay019,Bay2023,BJT019,CF2023,CF2023b,F2,FH2023,HMRZ021b}. The main purpose of this theory is to understand those properties shared by all homological shift ideals of a given monomial ideal. We call these properties, the \textit{homological shift properties}.

One of the driving motivation in this line of research is the \textit{Bandari-Bayati-Herzog conjecture} \cite{HMRZ021a} which asserts that the homological shift ideals of a polymatroidal are again polymatroidal. This conjecture is widely open. However, the conjecture was proved by Bayati for squarefree polymatroidal ideals \cite{Bay019}, by Herzog, Moradi, Rahimbeigi and Zhu for polymatroidal ideals that satisfy the strong exchange property \cite{HMRZ021a}, and by the author and Herzog for polymatroidal ideals generated in degree 2 \cite{FH2023}. Furthermore, it was shown by the author in \cite[Theorem 2.2]{F2} that $\HS_1(I)$ is always polymatroidal if $I$ is such, pointing towards the validity of the conjecture in general. This latter result was also recently recovered by Bayati in \cite[Corollary 2.2]{Bay2023}.

Another interesting conjecture about the homological shifts of powers of the cover ideals of Cohen--Macaulay very well-covered graphs was recently formulated in \cite{CF2023,CF2023b} and proved in some special cases, including bipartite and whisker graphs.

In the present paper, we illustrate and explain how to use the \textit{Macaulay2} \cite{GDS} package \texttt{HomologicalShiftIdeals}. In Section \ref{sec2-FHSPack} the mathematical background needed to develop some of the algorithms of the package is explained. In Section \ref{sec3-FHSPack}, some examples are presented, illustrating how to use the functions of the package.

\section{Mathematical background}\label{sec2-FHSPack}

Let $S=K[x_1,\dots,x_n]$ be the standard graded polynomial ring over a field $K$. We set ${\bf x^a}=x_1^{a_1}\cdots x_n^{a_n}$ for ${\bf a}=(a_1,\dots,a_n)\in\ZZ_{\ge0}^n$. The vector ${\bf a}$ is called the \textit{multidegree} of ${\bf x^a}$, whereas $\deg({\bf x^a})=a_1+a_2+\dots+a_n$ is its \textit{degree}.

Let $I$ be a monomial ideal of $S$. Since $I$ is multigraded, the minimal free resolution is multigraded as well, say
$$
\FF\ :\ \ \ \cdots\rightarrow F_i\rightarrow\cdots\rightarrow F_0\rightarrow I\rightarrow0,
$$
where $F_i=\bigoplus_{\bf a}S(-{\bf a})^{\beta_{i,{\bf a}}(I)}$ for all $i$, and where $\beta_{i,{\bf a}}(I)$ is the $(i,{\bf a})$th multigraded Betti number of $I$. The vectors ${\bf a}\in\ZZ^n_{\ge0}$ such that $\beta_{i,{\bf a}}(I)\ne0$ are called the $i$th \textit{multigraded shifts} of $I$. The \textit{projective dimension} of $I$ is defined as the integer $\pd(I)=\max\{i:\beta_i(I)\ne0\}$. Whereas, the (\textit{Castelnuovo--Mumford}) \textit{regularity} of $I$ is the integer $\reg(I)=\max\{\deg({\bf x^a})-i:\beta_{i,{\bf a}}(I)\ne0\}$.
\begin{Definition}
\rm The \textit{$i$th homological shift ideal} of $I$ is the monomial ideal
$$
\HS_i(I)\ =\ ({\bf x}^{{\bf a}}\ :\ \beta_{i,{\bf a}}(I)\ne0).
$$
\end{Definition}
Note that $\HS_0(I)=I$ and $\HS_i(I)=(0)$ for $i<0$ and $i>\pd(I)$.\medskip

The main purpose of the theory is to determine those properties enjoyed by all $\HS_j(I)$. We call these properties, the \textit{homological shift properties} of $I$.

Let $I\subset S$ be a monomial ideal, and let $G(I)$ be its unique minimal monomial generating set. The \textit{initial degree} of $I$ is $\indeg(I)=\min\{\deg(u):u\in G(I)\}$.
\begin{Definition}
	\rm Let $I\subset S$ be a monomial ideal, and let $G(I)=\{u_1,\dots,u_m\}$.
	\begin{enumerate}[label=(\alph*)]
		\item $I$ has a \textit{linear resolution} if $\indeg(I)=\reg(I)$.
		\item $I$ has \textit{linear quotients} if there exists an order $u_1,\dots,u_m$ of $G(I)$ such that $(u_1,\dots,u_{k-1}):u_k$ is generated by a subset of the variables for $k=2,\dots,m$. In this case, $u_1,\dots,u_m$ is called an \textit{admissible order} of $I$.
		\item Let $u={\bf x^a}=x_1^{a_1}\cdots x_n^{a_n}$. The \textit{$x_i$-degree} of $u$ is the integer $\deg_{x_i}(u)=a_i$. We say that $I$ is \textit{polymatroidal} if $I$ is equigenerated and the \textit{exchange property} holds: for all $u,v\in G(I)$ and all $i$ with $\deg_{x_i}(u)>\deg_{x_i}(v)$ there exists $j$ such that $\deg_{x_j}(u)<\deg_{x_j}(v)$ and $x_j(u/x_i)\in G(I)$.
	\end{enumerate}
    If (a), or (b), or (c), is an homological shift property, we say that $I$: has \textit{homological linear resolution}, respectively, \textit{homological linear quotients}, respectively, is \textit{homological polymatroidal}.
\end{Definition}

For an equigenerated monomial ideal $I\subset S$, the following hierarchy holds:
\begin{center}
	(c) $\Rightarrow$ (b) $\Rightarrow$ (a).
\end{center}

Before we proceed, we recall some other concepts. The \textit{support} of a monomial $u\in S$ is the set $\supp(u)=\{x_i:\deg_{x_i}(u)>0\}$. Let $I\subset S$ be a monomial ideal. The \textit{support} of $I$ is the set $\supp(I)=\bigcup_{u\in G(I)}\supp(u)$. We say that $I$ is \textit{fully supported} (in $S$) if $\supp(I)=\{x_1,\dots,x_n\}$. The \textit{bounding multidegree} of $I$ is the vector $\textbf{deg}(I)=(\deg_{x_1}(I),\dots,\deg_{x_n}(I))\in\ZZ^n$, defined by
$$
\deg_{x_i}(I)\ =\ \max_{u\in G(I)}\deg_{x_i}(u).
$$

Furthermore, the \textit{socle} $\soc(I)$ of a monomial ideal $I\subset S$ is the set of monomials of $(I:\mathfrak{m})\setminus I$, where $\mathfrak{m}=(x_1,\dots,x_n)$. In other words, $\soc(I)$ is the set of all monomials $v$ such that $v\notin I$ and $x_iv\in I$, for $i=1,\dots,n$.\medskip

The purpose of the package \texttt{HomologicalShiftIdeals} is to provide the tools to manipulate and calculate the homological shift ideals of a monomial ideal $I$ and to determine the homological shift properties of $I$. The next table collects the functions available in the package and their use. $I\subset S$ denotes a monomial ideal, ${\bf a}$ an integral vector, ${\bf x^a}$ a monomial, $i$ an integer and $L$ a list of monomials.
\small\begin{table}[H]
	\centering
	\begin{tabular}{ll}
		\rowcolor{black!20}\bottomrule[1.05pt]
		Functions&Description\\
		\toprule[1.05pt]
		\texttt{supportIdeal}$(I)$&Computes $\supp(I)$\\
		\texttt{isFullySupported}$(I)$&Checks whether $I$ is fully supported (in $S$)\\
		\texttt{toMonomial$(S,{\bf a})$}&Computes the monomial ${\bf x^a}$ if ${\bf a}\in\ZZ^n_{\ge0}$\\
		\texttt{toMultidegree}$({\bf x^a})$&Computes the multidegree ${\bf a}$ of ${\bf x^a}$\\
		\texttt{boundingMultidegree}$(I)$&Computes $\textbf{deg}(I)$\\
		\texttt{multigradedShifts}$(I,i)$&Computes the $i$th multigraded shifts of $I$\\
		\texttt{HS}$(I,i)$&Computes $\HS_i(I)$\\
		\texttt{socle}$(I)$&Computes $\soc(I)$\\
		\texttt{hasLinearResolution}$(I)$&Checks if $I$ has a linear resolution\\
		\texttt{hasHomologicalLinearResolution}$(I)$&Checks if $I$ has homological linear resolution\\
		\texttt{hasLinearQuotients}$(I)$&Checks if $I$ has linear quotients\\
		\texttt{hasHomologicalLinearQuotients}$(I)$&Checks if $I$ has homological linear quotients\\
		\texttt{admissibleOrder}$(I)$&Determines an admissible order of $I$\\
		\texttt{isAdmissibleOrder}$(I,L)$&Checks if $L$ is an admissible order of $I$\\
		\texttt{isPolymatroidal}$(I)$&Checks if $I$ is polymatroidal\\
		\texttt{isHomologicalPolymatroidal}$(I)$&Checks if $I$ is homological polymatroidal\\
		\bottomrule[1.05pt]
	\end{tabular}\medskip
	\caption{List of the functions of \texttt{HomologicalShiftIdeals}.}
\end{table}
\normalsize

In the remaining part of this section, we explain the theory behind some of the algorithms used in the package. Given $n\in\mathbb{N}$, we set $[n]=\{1,\dots,n\}$. For a nonempty subset $A$ of $[n]$, we set ${\bf x}_A=\prod_{i\in A}x_i$.

We start with the function \texttt{socle}.
\begin{Proposition}
	\textup{\cite[Proposition 1.13]{HMRZ021a}} Let $I\subset S$ be a monomial ideal. Then $$\soc(I)\ =\ \{{\bf x^a}\ :\ {\bf x^a}\in(I:\mm)\setminus I\}.$$
    In particular, $\beta_{n-1,{\bf a}}(I)\ne0$ if and only if ${\bf x^a}/{\bf x}_{[n]}\in\soc(I)$.
\end{Proposition}

The previous result justifies the following algorithm that calculates \texttt{socle}$(I)$.
\begin{enumerate}
	\item[] \begin{enumerate}
	\item[\textsc{Step 1:}] Compute $M=$ \texttt{multigradedShifts}$(I,n-1)$.
	\item[\textsc{Step 2:}] Compute $\soc(I)=\{w/(x_1x_2\cdots x_n):w\in M\}$.
	\end{enumerate}
\end{enumerate}\medskip

Next, we discuss the functions \texttt{hasLinearQuotients} and \texttt{admissibleOrder}. For their implementations, we have imported the packages \texttt{SimplicialComplexes}, and \texttt{SimplicialDecomposability}. These packages use the \textit{Alexander duality} theory in a fruitful way \cite{DC2010}.

A \textit{simplicial complex} on the \textit{vertex set} $[n]$ is a family of subsets of $[n]$ such that 
\begin{enumerate}
	\item[-] $\{i\}\in\Delta$ for all $i\in[n]$, and 
	\item[-] if $F\subseteq[n]$ and $G\subseteq F$, we have $G\in\Delta$. 
\end{enumerate}
The dimension of $\Delta$ is the number $d=\max\{|F|-1:F\in\Delta\}$. Any $F\in\Delta$ is called a \textit{face} and $|F|-1$ is the \textit{dimension} of $F$. A \textit{facet} of $\Delta$ is a maximal face with respect to the inclusion. The set of facets of $\Delta$ is denoted by $\mathcal{F}(\Delta)=\{F_1,\dots,F_m\}$. In this case we write $\Delta=\langle F_1,\dots,F_m\rangle$ and say that $F_1,\dots,F_m$ \textit{generates} $\Delta$. We say that $\Delta$ is \textit{pure} of dimension $d$ if all facets of $\Delta$ have dimension $d$. The \textit{Alexander dual} of $\Delta$ is the simplicial complex (see \cite[Lemma 1.5.2]{JT}) defined by
$$
\Delta^\vee\ =\ \{[n]\setminus F\ :\ F\notin\Delta\}.
$$

A monomial $u\in S$ is \textit{squarefree} if $\deg_{x_i}(u)\le 1$, for all $i\in[n]$. A monomial ideal $I\subset S$ is \textit{squarefree} if each $u\in G(I)$ is squarefree. It is well known that for any squarefree ideal $I\subset S$ there exists a unique simplicial complex $\Delta$ on $[n]$ such that $I=I_\Delta$, where $I_\Delta=({\bf x}_F:F\subseteq[n],F\notin\Delta)$ is the \textit{Stanley--Reisner ideal} of $\Delta$ \cite{JT}.

Now, we establish the connection between squarefree monomial ideals with linear quotients and the \textit{shellability} of simplicial complexes. Recall that $\Delta$ is \textit{shellable} if there exists an order $F_1,F_2,\dots,F_m$ of its facets $\mathcal{F}(\Delta)$ such that
$$
\langle F_1,\dots,F_{k-1}\rangle\cap\langle F_{k}\rangle
$$
is pure of dimension $\dim(F_k)-1$ for $k=2,\dots,m$. Any order of the facets of $\Delta$ satisfying the conditions above is called a \textit{shelling order} of $\Delta$. The following result shows that admissible orders and shelling orders are essentially the same thing.
\begin{Theorem}\label{Thm:AdmOrShellOrd}
	\textup{\cite[Proposition 8.2.5]{JT}} The following conditions are equivalent.
	\begin{enumerate}[label=\textup{(\alph*)}]
		\item $I_\Delta$ has linear quotients.
		\item The Alexander dual $\Delta^\vee$ of $\Delta$ is shellable.
	\end{enumerate}
    Furthermore, $F_1,\dots,F_m$ is a shelling order of the Alexander dual $\Delta^\vee$ of $\Delta$, if and only if, ${\bf x}_{[n]\setminus F_1},\dots,{\bf x}_{[n]\setminus F_m}$ is an admissible order of $I_\Delta$.
\end{Theorem}
The previous theorem provides an algorithm to determine an admissible order of a squarefree monomial ideal with linear quotients. In order to extend the above result to all (non necessarily squarefree) monomial ideals we use \textit{polarization}.

Let $u={\bf x^a}$ be a monomial of $S$. The \textit{polarization} of $u$ is the monomial
$$
	u^\wp\ =\ \prod_{i=1}^{n}(\prod_{j=1}^{a_i}x_{i,j})\ =\ \prod\limits_{\substack{i=1,\dots,n\\ a_i>0}}x_{i,1}x_{i,2}\cdots x_{i,a_i}.
$$

Let $S^\wp=K[x_{i,j}:i\in[n],j\in[\deg_{x_i}(I)]]$. The \textit{polarization} of $I$ is the monomial ideal $I^\wp\subset S^\wp$ with minimal generating set $G(I^\wp)=\{u^\wp:u\in G(I)\}$. 

The following lemma, taken from \cite[Lemma 4.10]{CF2023}, is pivotal.
\begin{Lemma}\label{Lem:IwpSets}
	Let $I\subset S$ be a monomial ideal with $G(I)=\{u_1,\dots,u_m\}$ and having linear quotients. Then, $u_1,\dots,u_m$ is an admissible order of $I$ if and only if $u_1^\wp,\dots,u_m^\wp$ is an admissible order of $I^\wp$. In particular, $I$ has linear quotients if and only if $I^\wp$ has linear quotients.
\end{Lemma}

Theorem \ref{Thm:AdmOrShellOrd} and Lemma \ref{Lem:IwpSets} justify the next algorithm for \texttt{hasLinearQuotients}, that determines whether a monomial ideal has linear quotients or not.
\begin{enumerate}
	\item[]
	\begin{enumerate}
	\item[\textsc{Step 1:}] Compute $I^\wp$.
	\item[\textsc{Step 2:}] Using \texttt{SimplicialComplexes} compute the Alexander dual $\Delta^\vee$, where $I^\wp=I_\Delta$.
	\item[\textsc{Step 3:}] Using \texttt{SimplicialDecomposability} determine if $\Delta^\vee$ is shellable. If the answer is yes then \texttt{hasLinearQuotients}$(I)=$ \texttt{true}, otherwise $=$ \texttt{false}.
	\end{enumerate}
\end{enumerate}\smallskip

For the function \texttt{admissibleOrder} as above we begin with \textsc{Step} 1, 2 and 3. If \texttt{hasLinearQuotients}$(I)=$ \texttt{false}, then $I$ does not have linear quotients. Otherwise, we complete our algorithm with the next two steps.
\begin{enumerate}
	\item[]
	\begin{enumerate}
		\item[\textsc{Step 4:}] Using \texttt{SimplicialDecomposability} compute a shelling order $F_1,\dots,F_m$ of the Alexander dual $\Delta^\vee$, where $I_\Delta=I^\wp$.
		\item[\textsc{Step 5:}] Determine the associated admissible order $u_1^\wp,\dots,u_m^\wp$ of $I^\wp$ and \textit{depolarize} it (by the substitutions $x_{i,j}\mapsto x_i$ for all $i$ and $j$) to obtain an admissible order of $I$.
	\end{enumerate}
\end{enumerate}
\section{Examples}\label{sec3-FHSPack}

In this final section, we present some examples to illustrate how to use the package.\smallskip

Let $I = (abd, abf, ace, adc, aef, bde, bcf, bce, cdf, def)\subset\QQ[a,\dots,f]$. Using the package we can check that $I$ has linear resolution but not linear quotients as follows.\medskip

\texttt{\phantom{i}i1: S = QQ[a..f];}

\texttt{\phantom{i}i2: I = ideal(a*b*d, a*b*f, a*c*e, a*d*c, a*e*f, b*d*e, b*c*f,}

\texttt{\phantom{i}\phantom{i2: }b*c*e, c*d*f, d*e*f);}

\texttt{\phantom{i}i3: loadPackage "HomologicalShiftIdeals"}

\texttt{\phantom{i}i4: hasLinearResolution I}

\texttt{\phantom{i}o4: true}

\texttt{\phantom{i}i5: hasLinearQuotients I}

\texttt{\phantom{i}o5: false}\medskip

In \cite[Theorem 1.4]{FH2023} we proved that for any monomial ideal generated in a single degree and having linear quotients, then $\HS_1(I)$ has linear quotients as well. However, the higher homological shift ideals of $I$ may fail to have linear quotients.\smallskip

Consider the ideal $J=(ab,ac,ad,de,df)$ of $S=\QQ[a,\dots,f]$. Then $J$ has linear resolution, indeed it has linear quotients, and $\HS_1(J)$ has linear quotients as well. However, $\HS_2(J)$ does not have linear quotients, not even linear resolution.\medskip

\texttt{\phantom{i}i6: J = ideal(a*b, a*c, a*d, d*e, d*f);}

\texttt{\phantom{i}i7: HS(J,0)==J}

\texttt{\phantom{i}o7: true}

\texttt{\phantom{i}i8: HS(J,1)}

\texttt{\phantom{i}o8: ideal(a*b*c, a*b*d, a*c*d, a*d*e, a*d*f, d*e*f)}

\texttt{\phantom{i}i9: hasLinearQuotients HS(J,1)}

\texttt{\phantom{i}o9: true}

\texttt{i10: HS(J,2)}

\texttt{o10: ideal(a*b*c*d, a*d*e*f)}

\texttt{i11: hasLinearResolution HS(J,2)}

\texttt{o11: false}\medskip

Consider the principal Borel ideal $I=B(x_2^2x_3)$ of $S=\QQ[x_1,x_2,x_3]$. Then $I=(x_1^3,x_1^2x_2,x_1^2x_3,x_1x_2^2,x_1x_2x_3,x_2^3,x_2^2x_3)$ has homological linear quotients, indeed $I$ is even homological polymatroidal, see \cite[Theorem 3.4]{BJT019}.\medskip

\verb|i12:  S = QQ[x_1..x_3];|

\verb|i13:  I = ideal(x_1^3, x_1^2*x_2, x_1^2*x_3, x_1*x_2^2, x_1*x_2*x_3, |

\phantom{i13: ,,,,}\verb|x_2^3, x_2^2*x_3);|

\verb|i14:  hasHomologicalLinearQuotients I|

\verb|o14:  true|

\verb|i15:  admissibleOrder HS(I,2)|

\verb|o15:  {x_1^3*x_2*x_3, x_1^2*x_2^2*x_3, x_1*x_2^3*x_3}|

\verb|i16:  socle I|

\verb|o16:  {x_1^2, x_1*x_2, x_2^2}|

\verb|i17:  isHomologicalPolymatroidal I|

\verb|o17:  true|

\end{document}